\documentclass[a4paper,reqno,10pt]{amsart}

\usepackage{comment} 

\usepackage{hyphenat}

\usepackage{manfnt} 
\reversemarginpar 

\usepackage[svgnames]{xcolor}

\usepackage{url}

\usepackage{footmisc}

\usepackage{amsmath}
\usepackage{amsthm}
\usepackage{amssymb}

\allowdisplaybreaks 

\usepackage{amscd}

\usepackage[T1]{fontenc}
\usepackage{lmodern}

\usepackage{bbm}

\usepackage[scaled=.98]{BOONDOX-uprscr}

\usepackage[OMLmathsfit]{isomath}

\usepackage{mathabx}

\usepackage{leftindex}

\usepackage{enumitem}
\setlist{noitemsep}

\usepackage{stackengine}

\swapnumbers 

\newtheoremstyle{exercise}
  {3pt} 
  {3pt} 
  {\small\rmfamily} 
  {
} 
  {\rmfamily\scshape} 
  {.} 
  {.5em} 
  {} 

\newtheoremstyle{newplain}
  {5pt}
  {5pt}
  {\itshape}
  {}
  {\rmfamily\scshape}
  { ---}
  {.5em}
  {}

\newtheoremstyle{newremark}
  {5pt}
  {5pt}
  {\rmfamily}
  {}
  {\rmfamily\scshape}
  { ---}
  {.5em}
  {}

\theoremstyle{newplain}
\newtheorem*{Theorem*}{Theorem} 

\theoremstyle{newplain}
\newtheorem{Theorem}{Theorem}

\newtheorem{Proposition}[Theorem]{Proposition}
\newtheorem{Conjecture}[Theorem]{Conjecture}
\newtheorem{Definition}[Theorem]{Definition}

\theoremstyle{newremark}
\newtheorem{Empty}[Theorem]{}
\newtheorem{Remark}[Theorem]{Remark}

\newtheorem{Claim}[Theorem]{Claim}

\theoremstyle{exercise}

\numberwithin{Theorem}{section}

\newcommand{\R}{\mathbb{R}}

\newcommand{\calC}{\mathscr{C}}

\newcommand{\calF}{\mathscr{F}}

\newcommand{\calJ}{\mathscr{J}}

\newcommand{\calL}{\mathscr{L}}

\newcommand{\calV}{\mathscr{V}}

\newcommand{\bdelta}{\boldsymbol{\delta}}

\newcommand{\boldeta}{\boldsymbol{\eta}}

\newcommand{\bL}{\mathbf{L}}

\newcommand{\bT}{\mathbf{T}}

\newcommand{\fb}{\pmb{f}}  
\newcommand{\bg}{\pmb{g}}

\DeclareMathOperator{\rmbody}{\mathrm{body}}          
\DeclareMathOperator{\rmdiam}{\mathrm{diam}}          

\DeclareMathOperator{\rmdist}{\mathrm{dist}}          
\DeclareMathOperator{\rmint}{\mathrm{int}}            

\DeclareMathOperator{\rmsign}{\mathrm{sign}}          

\newcommand{\blseg}{\pmb{[}}
\newcommand{\brseg}{\pmb{]}}

\newcommand{\lno}{\left\bracevert \mkern -5mu}
\newcommand{\rno}{\mkern -5mu \right\bracevert}

\def\XXint#1#2#3{{%
\setbox0=\hbox{$#1{#2#3}{\int}$}
\vcenter{\hbox{$#2#3$}}\kern-.5\wd0}}

\newcommand{\KH}{K\!H}
\newcommand{\bKH}{\mathbf{\KH}}

\newcommand{\cqfd} {
\renewcommand{\qedsymbol}{$\blacksquare$}
\qed
\renewcommand{\qedsymbol}{$\square$} }

\newcommand{\veps}{\varepsilon}

\newcommand{\la}{\langle}
\newcommand{\ra}{\rangle}

\newcommand{\ts}{\mkern 3mu} 

\newcommand{\ie}{{\it i.e.}\ }
\newcommand{\eg}{{\it e.g.}\ }

\renewcommand{\em}{\bf}

\renewcommand{\leq}{\leqslant}
\renewcommand{\geq}{\geqslant}

\renewcommand{\subset}{\subseteq}

\newlength{\drop}

\newcommand*{\plogo}{\fbox{$\mathcal{HSH}$}}
\usepackage{trajan}

\newcommand*{\titleAT}{\begingroup
\drop=0.1\textheight
\vspace*{\drop}
\rule{0.93\textwidth}{1pt}\par 
\vspace{2pt}\vspace{-\baselineskip}
\rule{0.93\textwidth}{0.4pt}\par 
\vspace{0.5\drop}
\centering

\textcolor{Blue}{
{\trjnfamily
	\Huge MORE SO D'ANALYSE}\\[5\baselineskip]
{\trjnfamily
	\Large VOLUME I}\\[1.8\baselineskip]
{\trjnfamily
	\Huge MEASURE THEORY}\\[5\baselineskip]
}\par 
\vspace{0.25\drop}
\rule{0.3\textwidth}{0.4pt}\par 
\vspace{\drop} 

{\Large \scshape Thierry De Pauw}\par 
\vfill 
{\large \textcolor{Blue}{\plogo}}\\[0.5\baselineskip]

{\large \scshape Home Sweet Home}\par 
\vspace*{\drop}
\endgroup}

\begin{document}

\title{Banach-Lamperti for Kurzweil-Henstock}

\author[Th. De Pauw]{Thierry De Pauw}

\address{Institute for Theoretical Sciences / School of Science, Westlake University\\
No. 600, Dunyu Road, Xihu District, Hangzhou, Zhejiang, 310030, China}

\email{thierry.depauw@westlake.edu.cn}




\begin{abstract}
We identify isometric isomorphisms of the space of Kurzweil-Hen\-stock integrable functions as bi-absolutely-continuous changes of variable. 
\end{abstract}

\maketitle



\tableofcontents


\section{Introduction}

By a {\em cell} we mean a non-degenerate compact interval in $\R$.
If $I \subset \R$ is a cell then $\KH(I)$ denotes the space of Kurzweil-Henstock integrable functions on $I$ and $\bKH(I)$ is the corresponding space of equivalence classes almost everywhere of members of $\KH(I)$.
We denote by $f$ typical members of $\KH(I)$, by $\fb$ typical members of $\bKH(I)$, and by $\blseg f \brseg$ the equivalence class of $f$.
Thus, $\fb = \blseg f \brseg$ whenever $f \in \fb$.
We equip $\bKH(I)$ with Alexiewicz's norm (see \cite{ALE.48}) defined by
\begin{equation*}
\| \fb \|_A = \max \left\{ \left| \int_{[\min I,x]} f \right| : \min I < x \leq \max I \right\}
\end{equation*}
whenever $f \in \fb \in \bKH(I)$.
We study the isometries of this normed space.
Specifically, we prove the following.

\begin{Theorem*}
\label{main.thm}
Let $I$ and $\tilde{I}$ be two cells in $\R$ and $\bT : \bKH(\tilde{I}) \to \bKH(I)$ a linear operator.
The following are equivalent.
\begin{enumerate}
\item[(A)] $\bT$ is a surjective isometry.
\item[(B)] There exist $\sigma \in \{-1,1\}$ and an increasing homeomorphism $\phi : I \to \tilde{I}$ such that both $\phi$ and $\phi^{-1}$ are absolutely continuous and 
\begin{equation*}
\bT(\fb) = \blseg \sigma \cdot (f \circ \phi) \cdot \phi' \brseg ,
\end{equation*}
for all $f \in \fb \in \bKH(I)$.
\end{enumerate}
\end{Theorem*}
\par 
The proof that (B) implies (A) in section \ref{sec.cv} boils down to establishing the proper change of variable theorem in this context (see also \cite[Theorem 1]{LEA.03} and \cite[Theorem 7.3]{BEN.10}).
The proof of the reverse implication in section \ref{sec.iso} relies on the Banach-Stone theorem at the level of indefinite integrals together with the appropriate absolute continuity property of indefinite Kurzweil-Henstock integrals.
Recall that $\bKH(I)[\|\cdot\|_A]$ is not a Banach space (though, it is barrelled but this will not be used here), therefore, a linear isometric operator may have a dense image without being surjective.
\par 
This is in contrast with Banach-Lamperti's theorem (see \cite[Ch. XI, \S5, Espaces $L^{(p)}$, n°2, I]{BANACH} and \cite{LAM.58}) asserting that the surjective linear isometries $\bT : \bL_1(\tilde{I}) \to \bL_1(I)$ correspond to pairs $(\phi,h)$, where $\phi : I \to \tilde{I}$ is a bi-Lebesgue-measurable bijection and $h : I \to \R$ is Lebesgue-measurable such that $\mu = \phi_* \lambda << \lambda$ and $|h| = \frac{d\mu}{d\lambda} > 0$, by means of the formula $\bT(\blseg f \brseg) = \blseg (f \circ \phi) \cdot h \brseg$.
Here, $\lambda$ is the Lebesgue measure and $\frac{d\mu}{d\lambda}$ is a Radon-Nikod\'ym derivative.
Thus, the case of Kurzweil-Henstock integration corresponds, so to speak, to more regular changes of variable than the case of Lebesgue integration.
Observe, however, that from a functional-analytic point of view the two problems have little in common, as $\bL_1(I)$ is an $\calL_1$-space whereas $\bKH(I)$ is an $\calL_\infty$-space (see \eg \cite[Appendix F]{BENYAMINI.LINDENSTRAUSS} for the relevant definition).
\par 
In an attempt to reach a large audience, we do not assume any prior knowledge of Kurzweil-Henstock integrability (the interested read may consult \eg \cite{BARTLE.01,LEE.VYBORNY.00,MAWHIN,PFEFFER.93,SWARTZ}).
In particular, none of the results in the next section appear to be new but we have gathered and proved them in a fashion most appropriate for our applications.
Of course, we assume familiarity with the theory of Lebesgue measure and integration in the real line.  
\par 
It is a pleasure to acknowledge helpful discussions with Ph. Bouafia during the preparation of this paper.

\section{Preliminaries}

\begin{Empty}[Vocabulary and notations]
\label{2.1}
Two cells $I_1$ and $I_2$ are {\em non-overlapping} if $I_1 \cap I_2$ is either empty or a singleton.
A {\em $P$-family} $\Pi$ is a finite set of pairs $(J,x)$, where $J$ is a cell and $x \in J$, such that $J_1$ and $J_2$ are non-overlapping whenever $(J_1,x_1), (J_2,x_2) \in \Pi$ and $J_1 \neq J_2$.
The {\em body} of a $P$-family $\Pi$ is the set $\rmbody(\Pi) = \cup \{ J : (J,x) \in \Pi \}$.
We say that $\Pi$ is {\em in a cell} $I$ if $\rmbody(\Pi) \subset I$.
A $P$-family $\Pi$ is a {\em $P$-division} of a cell $I$ if $\rmbody(\Pi) = I$.
The {\em Riemann sum} associated with a $P$-family $\Pi$ is
\begin{equation*}
S(\Pi,f) = \sum_{(J,x) \in \Pi} f(x) \cdot \lno J \rno,
\end{equation*}
where $\lno J \rno$ is the length of $J$, \ie $\max J - \min J$.
A {\em gauge} on a set $E \subset \R$ is a positive function $\bdelta : E \to \R^+ \setminus \{0\}$.
If $\bdelta$ is a gauge on a cell $I$ and $\Pi$ is a $P$-family in $I$ then we say that $\Pi$ is {\em $\bdelta$-fine} whenever $J \subset [x-\bdelta(x),x+\bdelta(x)]$ for all $(J,x) \in \Pi$.
\end{Empty}

\begin{Empty}[$KH$-integral]
\label{2.2}
Let $I$ be a cell and $f : I \to \R$ a function.
We say that $f$ is {\em Kurzweil-Henstock integrable} in $I$ if there exists a real number $r$ satisfying the following property.
For every $\veps > 0$ there exists a gauge $\bdelta$ on $I$ such that $|S(\Pi,f) - r| < \veps$ for all $\bdelta$-fine $P$-divisions $\Pi$ of $I$.
In that case, $r$ is unique and denoted $\int_I f$.
The set $KH(I)$ consisting of all Kurzweil-Henstock integrable functions defined on $I$ is a vector space and $KH(I) \to \R : f \mapsto \int_I f$ is linear.
If $f \in KH(I)$ and $\hat{I} \subset I$ is a cell then the restriction $f|_{\hat{I}} \in KH(\hat{I})$.
\par 
We now state what is known as the {\em Saks-Henstock} theorem (see \eg \cite[5.4]{BARTLE.01}, \cite[CH.10, \S8]{MAWHIN}, or \cite[Ch.3, lemma 1]{SWARTZ}).
If $f \in KH(I)$ and $\veps > 0$ then there exists a gauge $\bdelta$ on $I$ such that 
\begin{equation*}
\sum_{(J,x) \in \Pi} \left| f(x) \cdot \lno J \rno - \int_J f|_J \right| < \veps
\end{equation*}
for all $\bdelta$-fine $P$-family $\Pi$ in $I$.
\end{Empty}

\begin{Empty}[Indefinite integral]
\label{2.3}
If $f \in KH(I)$ we define a function $\int^\bullet f : I \to \R$ by means of the formula
\begin{equation*}
\left( \int^\bullet f \right)(x) = \begin{cases}
0 & \text{if } x = \min I \\
\int_{I_x} f|_{I_x} & \text{otherwise},
\end{cases}
\end{equation*}
where $I_x = [\min I,x]$.
The function $\int^\bullet f$ is continuous, according to Hake's theorem (see \eg \cite[5.6]{BARTLE.01}, \cite[Ch.11 \S1]{MAWHIN}, or \cite[Ch.3, corollary 2]{SWARTZ}), and we call it the {\em indefinite integral} of $f$.
\end{Empty}

\begin{Empty}[Additive functions]
\label{2.4}
Let $\calC(I)$ be the set of cells contained in a cell $I$.
We say that $F : \calC(I) \to \R$ is {\em additive} whenever the following holds.
For all positive integer $N$ and $J,J_1,\ldots,J_N \in \calC(I)$ if $J = \cup_{n=1}^N J_n$ and $J_1,\ldots,J_N$ are non-overlapping then $F(J) = \sum_{n=1}^N F(J_n)$.
\par 
If $F : I \to \R$ then $F^\times : \calC(I) \to \R$ defined by $F^\times(J) = F(\max J) - F(\min J)$ is additive.
If $F : \calC(I) \to \R$ is additive then we let $\leftindex^\times {F} : I \to \R$ be defined by $\leftindex^\times {F}(\min I)=0$ and $\leftindex[I]^\times {F}(x) = F([\min I,x])$ whenever $\min I < x \leq \max I$.
Note that\footnote{In this sentence, the first two occurrences of the symbol $F$ and the last three refer to objects of different nature.} $(\leftindex^\times F)^\times = F$ and that $\leftindex[I]^\times{(F^\times)}=F$ if and only if $F(\min I)=0$.
This establishes a bijection from $\R^I \cap \{ F : F(\min I) = 0 \}$ to the set of additive functions $\calC(I) \to \R$.
\par 
The ``indefinite integral'' $\left(\int^\bullet f\right)^\times$ of $f \in KH(I)$ is additive.
The less pedantic choice not to distinguish between $F$ and $F^\times$ leads to oddities such as ``$F$ is increasing if and only if $F$ is positive''.
\end{Empty}

\begin{Empty}[Lebesgue measure]
\label{2.6}
If $A \subset \R$ we let $\lno A \rno$ be its Lebesgue outer measure.
When we say that $A$ is {\em negligible} it always refer to Lebesgue measure, \ie $\lno A \rno = 0$.
The phrase ``almost everywhere'' always refers to Lebesgue-negligibility.
\end{Empty}

\begin{Empty}[The condition $AC_*$]
\label{2.5}
A $P$-family $\Pi$ in $I$ is said to be {\em tagged} in a set $E \subset I$ if $x \in E$ for all $(J,x) \in \Pi$.
We say that an additive function $F : \calC(I) \to \R$ is $\boldsymbol{AC_*}$ if it satisfies the following property.
For every negligible set $N \subset I$ and every $\veps > 0$ there exists a gauge $\bdelta$ on $N$ such that
\begin{equation*}
\sum_{(J,x) \in \Pi} |F(J)| < \veps
\end{equation*}
for all $\bdelta$-fine $P$-families in $I$ tagged in $N$.
\end{Empty}

\begin{Theorem}
\label{2.7}
Let $F : \calC(I) \to \R$ be an additive function.
The following are equivalent.
\begin{enumerate}
\item[(A)] There exists $f \in \KH(I)$ such that $F = \left(\int^\bullet f\right)^\times$.
\item[(B)] $F$ is $AC_*$ and $\leftindex^\times F$ is differentiable almost everywhere. 
\end{enumerate}
Furthermore, in case (B) holds, $(\leftindex^\times F)' \in \KH(I)$ and $\leftindex^\times F = \int^\bullet (\leftindex^\times F)'$.
\end{Theorem}

\begin{proof}[Proof that $(A) \Rightarrow (B)$]
{\bf (i)}
Let $N \subset I$ be negligible and $\veps > 0$.
For each natural integer $j$ define $N_j = N \cap \{ j \leq |f| < j+1 \}$ and choose an open set $U_j$ containing $N_j$ and such that $\lno U_j \rno < \frac{\veps}{(j+1) \cdot 2^{j+2}}$.
Let $\boldeta$ be a gauge on $I$ associated to $f$ and $\frac{\veps}{2}$ by means of the Saks-Henstock theorem.
Define a gauge $\bdelta$ on $N$ by setting $\bdelta|_{N_j} = \min \{ \frac{1}{2} \rmdist(\cdot , \R \setminus U_j),\boldeta \}$ for all $j$.
\par 
Let $\Pi$ be a $\bdelta$-fine $P$-family in $I$ tagged in $N$.
Since $\Pi$ is $\boldeta$-fine, it follows from the Saks-Henstock theorem that
\begin{equation*}
\sum_{(J,x) \in \Pi} |F(J)| < \frac{\veps}{2} + \sum_{(J,x) \in \Pi} |f(x)|\cdot \lno J \rno.
\end{equation*}
Moreover,
\begin{multline*}
\sum_{(J,x) \in \Pi} |f(x)|\cdot \lno J \rno = \sum_{j = 0}^\infty \sum_{\substack{(J,x) \in \Pi \\ x \in N_j}} |f(x)|\cdot \lno J \rno \leq \sum_{j = 0}^\infty (j+1) \cdot \sum_{\substack{(J,x) \in \Pi \\ x \in N_j}} \lno J \rno  \\ \leq \sum_{j=0}^\infty (j+1) \cdot \lno U_j \rno < \frac{\veps}{2},
\end{multline*}
since $\rmbody(\Pi_j) \subset U_j$, where $\Pi_j = \Pi \cap \{ (J,x) : x \in N_j \}$, for all $j$, by the choice of $\bdelta|_{N_j}$.
This completes the proof that $F$ is $AC_*$.
\par 
{\bf (ii)}
Let $G = \rmint(I) \cap \{ x : F \text{ is differentiable at $x$ and } F'(x) = f(x) \}$ and $B = \rmint(I) \setminus G$.
We shall show that $B$ is negligible.
To this end, we define 
\begin{equation*}
B_k = \rmint(I) \cap \left\{ x : \limsup_{\substack{h \to 0 \\ h \neq 0}} \left| \frac{F(x+h)-F(x)}{h} - f(x) \right| > \frac{1}{k} \right\} 
\end{equation*}
for each positive integer $k$ and we notice that $B = \cup_k B_k$.
Thus, it suffices to prove that each $B_k$ is negligible.
\par 
Fix $k$.
Given $\veps > 0$, choose a gauge $\bdelta$ on $I$ associated with $f$ and $\veps$ in Saks-Henstock's theorem.
Define 
\begin{multline}
\label{eq.1}
\calV = \calC(I) \cap \Bigg\{ J : \big|F(J) - f(\min J)\cdot \lno J \rno \big| > \frac{\lno J \rno}{k} \text{ and } \lno J \rno \leq \bdelta(\min J) \textbf{ or } \\
\big|F(J) - f(\max J)\cdot \lno J \rno \big| > \frac{\lno J \rno}{k} \text{ and } \lno J \rno \leq \bdelta(\max J) \Bigg\}.
\end{multline}
Observe that for all $x \in B_k$ there exists $J \in \calV$ of which $x$ is an endpoint. 
By the Vitali covering theorem, there are finitely many pairwise disjoint members of $\calV$, say $J_1,\ldots,J_N$, such that $\lno B_k \setminus \cup_{n=1}^N J_n \rno < \veps$.
For each $n$, let $x_n = \min J_n$ or $x_n = \max J_n$ according to whether the first or second condition in \eqref{eq.1} witnesses membership of $J_n$ to $\calV$.
Thus, the $P$-family $\{ (J_n,x_n) : n=1,\ldots,N \}$ is $\bdelta$-fine and we infer from the Saks-Henstock theorem that
\begin{equation*}
\sum_{n=1}^N \lno J_n \rno \leq k \cdot \sum_{n=1}^N \big| F(J_n) - f(x_n) \cdot \lno J_n \rno \big| < k \cdot \veps.
\end{equation*}
It ensues that $\lno B_k \rno \leq \lno B_k \setminus \cup_{n=1}^N J_n \rno + \lno \cup_{n=1}^N J_n \rno < (1 + k) \cdot \veps$.
Since $\veps > 0$ is arbitrary, the proof is complete.
\end{proof}

\begin{Remark}
\label{2.8}
The following comments are in order.
\begin{enumerate}
\item[(1)] In the above proof that (A) implies (B), we have shown, in fact, that if $f \in KH(I)$ then $\left( \int^\bullet f \right)'(x) = f(x)$ for almost every $x \in I$.
\item[(2)] If $f \in KH(I)$ and $\left( \int^\bullet f \right)(x)=0$ for all $x \in I$ then $f(x)=0$ for almost every $x \in I$.
Indeed, under this assumption the function $x \mapsto \int_{\min I}^x f$ is differentiable everywhere and its derivative vanishes identically.
The conclusion then ensues from (1).
\item[(3)] Arguing as in {\bf (i)} above, one shows that if $f : I \to \R$ and $\{f \neq 0 \}$ is negligible then $f \in KH(I)$ and $\int_I f = 0$.
In particular, if $f,g : I \to \R$ and $\{ f \neq g \}$ is negligible then $f \in KH(I)$ if and only if $g \in KH(I)$, in which case $\int_I f = \int_I g$.
This clarifies the last sentence in the statement of \ref{2.7} as well as the statement of our main theorem in the introduction (as $\left(\leftindex^\times F\right)'$, respectively $\phi'$, is defined merely almost everywhere).
\end{enumerate}
\end{Remark}

\begin{proof}[Proof that $(B) \Rightarrow (A)$]
Let $N \subset I$ consist of the endpoints of $I$ and of all points of $\rmint(I)$ at which $F$ is not differentiable.
Thus, $N$ is negligible.
Define $f : I \to \R$ by declaring that $f(x) = (\leftindex^\times F)'(x)$ if $x \in I \setminus N$ and $f(x)=0$ if $x \in N$.
We will show that $f \in KH(I)$ and $\int_I f = F(I)$.
Repeating this for all cells contained in $I$ will complete the proof.
\par 
Fix $\veps > 0$ and $x \in I \setminus N$.
Since $\leftindex^\times F$ is differentiable at $x$ and its derivative is $f(x)$, one checks that there exists $\bdelta_{\mathrm{diff}}(x) > 0$ such that, for all $J \in \calC(I)$, if\footnote{The condition $x \in J$ -- in the definition of Kurzweil-Henstock integrability -- is essential precisely here.} $x \in J$ and $\lno J \rno < \bdelta_{\mathrm{diff}}(x)$ then $\big| F(J) - f(x) \cdot \lno J \rno \big| < \veps \cdot \lno J \rno$.
Furthermore, there exists a gauge $\bdelta_{\mathrm{ac}}$ on $N$ associated with $F$ and $\veps$ is the definition of $F$ being $AC_*$.
We now define a gauge $\bdelta$ on $I$ as follows:
\begin{equation*}
\bdelta(x) = \begin{cases}
\bdelta_{\mathrm{ac}}(x)  & \text{if } x \in N \\
 \bdelta_{\mathrm{diff}}(x) & \text{if } x \in I \setminus N .
\end{cases}
\end{equation*}
\par 
If $\Pi$ is any $P$-division of $I$ then
\begin{gather*}
F(I) = \sum_{(J,x) \in \Pi} F(J) = \sum_{\substack{(J,x) \in \Pi \\ x \not\in N}} F(J) +  \sum_{\substack{(J,x) \in \Pi \\ x \in N}} F(J) ,\\
S(\Pi,f) = \sum_{(J,x) \in \Pi} f(x) \cdot \lno J \rno = \sum_{\substack{(J,x) \in \Pi \\ x \not\in N}} f(x) \cdot \lno J \rno.
\end{gather*}
Therefore,
\begin{equation*}
\big| F(I) - S(\Pi,f) \big| \leq \sum_{\substack{(J,x) \in \Pi \\ x \not\in N}} \big| F(J) - f(x) \cdot \lno J \rno \big| + \sum_{\substack{(J,x) \in \Pi \\ x \in N}} |F(J)|.
\end{equation*}
If, moreover, $\Pi$ is $\bdelta$-fine then
\begin{equation*}
\big| F(I) - S(\Pi,f) \big| <  \left( \sum_{\substack{(J,x) \in \Pi \\ x \not\in N}} \veps \cdot \lno J \rno \right) + \veps \leq \left( \lno I \rno + 1\right) \cdot \veps.
\end{equation*}
\end{proof}

\begin{Empty}[The condition $AC$]
\label{2.9}
We say that an additive function $F : \calC(I) \to \R$ is $\boldsymbol{AC}$ whenever the following condition is satisfied.
For every $\veps > 0$ there exists $\delta > 0$ such that for all positive integers $N$ and all $J,J_1,\ldots,J_N \in \calC(I)$ if $J = \cup_{n=1}^N J_n$ and $J_1,\ldots,J_N$ are pairwise non-overlapping and $\sum_{n=1}^N \lno J_n \rno < \delta$ then $\sum_{n=1}^N |F(J_n)| < \veps$.
This is the classical notion pertaining to measure theory.
\begin{enumerate}
\item[(A)] {\it If $F$ is $AC$ then it is $AC_*$.}
\end{enumerate}
\par 
{\it Proof.}
Fix $\veps > 0$ and let $\delta > 0$ be as in the definition above.
Let $N \subset I$ be negligible and choose an open set $U \subset \R$ containing $N$ such that $\lno U \rno < \delta$.
Define a gauge $\bdelta$ on $N$ to coincide with $\frac{1}{2}\rmdist(\cdot , \R \setminus U)$.
If $\Pi$ is a $\bdelta$-fine $P$-family anchored in $N$ then $\rmbody(\Pi) \subset U$.
Therefore, $\sum_{(J,x) \in \Pi} \lno J \rno \leq \lno U \rno < \delta$ and, in turn, $\sum_{(J,x) \in \Pi} |F(J)| < \veps$.
\cqfd
\begin{enumerate}
\item[(B)] {\it The converse of (A) fails.}
\end{enumerate}
\par 
{\it Proof.}
This classical fact is witnessed for instance by $F^\times$ on $I = [0,1]$, where $F(x) = x^2 \sin \left( \frac{1}{x^2} \right)$ if $0 < x \leq 1$ and $F(0)=0$.\cqfd
\end{Empty}

The following is the 1-dimensional version of the Besicovitch covering theorem.
For this and its higher dimensional analogues, see \eg \cite[Ch. I \S1]{GUZ.75} or \cite[5.2.4]{PFEFFER.93} (the latter gives a proof leading to 14 instead of 5 subcollections but none of this matters for our subsequent application). 
By a {\em disjointed} collection of sets we mean that any two distinct members of it are disjoint.

\begin{Theorem}[A covering theorem]
\label{2.10}
Let $X \subset \R$ be bounded and let $\bdelta$ be a bounded gauge on $X$.
Abbreviate $J(x) = [x-\bdelta(x),x+\bdelta(x)]$.
Then there exist five disjointed collections $\calJ^*_i \subset \{ J(x) : x \in X \}$, $i=1,\ldots,5$, such that $X \subset \cup_{i=1}^5 \left( \cup \calJ^*_i \right)$.
\end{Theorem}

\begin{Empty}[Luzin's condition (N)]
\label{2.11}
We say that a function $g : I \to \R$ satisfies {\em Luzin's condition $(N)$} if $\lno g(Z) \rno = 0$ for all $Z \subset I$ such that $\lno Z \rno = 0$, \ie $g$ maps negligible sets to negligible sets.
\end{Empty}

The following is a special case of the Banach-Zarecki theorem (see \eg \cite[Ch. 7 \S3]{BRU.BRU.THO}).

\begin{Theorem}
\label{2.12}
Let $F : \calC(I) \to \R$ be additive and non-negative, \ie $F(J) \geq 0$ for all $J \in \calC(I)$.
The following are equivalent.
\begin{enumerate}
\item[(A)] $F$ is $AC_*$.
\item[(B)] $\leftindex^\times F$ is continuous and satisfies Luzin's condition $(N)$.
\item[(C)] $F$ is $AC$.
\end{enumerate}
\end{Theorem}

\begin{proof}
{\bf (i)}
We make a preliminary comment about how we shall use the hypothesis that $F$ be non-negative.
This is equivalent to $\leftindex^\times F$ being non-decreasing. 
Thus, for any cell $J = [a,b] \subset I$ we have $ \leftindex^\times F(J) =  \leftindex^\times F([a,b]) = \left[  \leftindex^\times F(a),  \leftindex^\times F(b)\right]$ so that
\begin{equation}
\label{eq.2}
\lno  \leftindex^\times F(J) \rno =  \leftindex^\times F(b) -  \leftindex^\times F(a) = F([a,b]) = F(J).
\end{equation}
Similarly, if $\hat{J} \subset \R$ is a cell then $J := \left( \leftindex^\times F\right)^{-1}(\hat{J}) \subset I$ is a cell as well and, upon writing $J = [a,b]$, we have
\begin{equation}
\label{eq.3}
F(J) = \leftindex^\times F(b) - \leftindex^\times F(a) \leq \lno \hat{J} \rno,
\end{equation}
Since $\leftindex^\times F(a), \leftindex^\times F(b) \in \hat{J}$.
\par 
{\bf (ii)} 
Here, we prove that (A) implies (B).
The continuity of $\leftindex^\times F$ readily follows from the definition of $F$ being $AC_*$.
Thus, given a negligible set $Z \subset I$ we ought to show that $\leftindex^\times F(Z)$ is negligible as well.
Given $\veps > 0$, the $AC_*$ property of $F$ implies the existence of a gauge $\bdelta$ on $Z$ such that $\sum_{(J,x) \in \Pi} |F(J)| < \veps$ whenever $\Pi$ is a $\bdelta$-fine $P$-family anchored in $Z$.
The same holds when $\bdelta$ is replaced with $\min\{\bdelta,\rmdiam I\}$.
Accordingly, there is no restriction to assume that $\bdelta$ is bounded.
Since 
\begin{equation*}
Z \subset \bigcup_{x \in Z}  [x-\bdelta(x),x+\bdelta(x)],
\end{equation*} 
in view of theorem \ref{2.10}, we infer the existence of five disjointed collections of cells $\calJ_i^*$, $i=1,\ldots,5$, each member of which is of the type $J_x = I \cap [x-\bdelta(x),x+\bdelta(x)]$ for some $x \in Z$, such that $Z \subset \cup_{i=1}^5 \left( \cup \calJ^*_i \right)$.
Given $i=1,\ldots,5$ and $\calF \subset \calJ_i^*$ a finite subcollection, we define $\Pi_\calF = \{ (J,x(J)) : J \in \calF \}$, where $x(J)$ denotes the center of $J$.
Note that $\Pi_\calF$ is a $\bdelta$-fine $P$-family anchored in $Z$, thus, $\sum_{J \in \calF} |F(J)| = \sum_{(J,x) \in \Pi_\calF} |F(J)| < \veps$.
As $\calF$ is arbitrary, we infer that $\sum_{J \in \calJ^*_i} |F(J)| \leq \veps$, $i=1,\ldots,5$.
\par 
Now, 
\begin{equation*}
\leftindex^\times F(Z) \subset \leftindex^\times F \left( \bigcup_{i=1}^5 \bigcup \calJ^*_i\right)
\subset \bigcup_{i=1}^5 \bigcup_{J \in \calJ^*_i} \leftindex^\times F(J)
\end{equation*}
and, in view of {\bf (i)} (specifically \eqref{eq.2}), we conclude that
\begin{equation}
\lno \leftindex^\times F(Z) \rno \leq \sum_{i=1}^5 \sum_{J \in \calJ_i^*} \lno \leftindex^\times F(J) \rno = \sum_{i=1}^5 \sum_{J \in \calJ_i^*} | F(J)|  \leq 5 \cdot \veps.
\end{equation}
Since $\veps$ is arbitrary, $\lno \leftindex^\times F(Z) \rno =0$.
\par 
{\bf (iii)}
We turn to showing that (B) implies (C).
We briefly recall the construction of the Lebesgue-Stieltjes measure $\mu$ associated with $\leftindex^\times F$.
For an arbitrary $A \subset I$ we put
\begin{equation*}
\mu^*(A) = \inf \left\{ \sum_{k=1}^\infty F(J_k) : \la J_k \ra_k \text{ is a sequence in } \calC(I) \text{ and } A \subset \bigcup_{k=1}^\infty J_k \right\}.
\end{equation*}
Borel subsets of $I$ are $\mu^*$-measurable in the sense of Carathéodory and $\mu^*(J) = F(J)$ for all $J \in \calC(I)$ (see \eg \cite[Theorem 3.20 conclusions 1 and 4]{BRU.BRU.THO}).
In particular, the restriction $\mu$ of $\mu^*$ to the set of Borel subsets of $I$ is a finite Borel measure on $I$.
Notice also that $\mu(\{x\})=0$ for all $x \in I$, since $\leftindex^\times F$ is continuous.
In particular, if $J_1$ and $J_2$ are non-overlapping cells contained in $I$ then $\mu(J_1 \cup J_2) = \mu(J_1) + \mu(J_2)$.
\par 
We now show that $\mu$ is absolutely continuous with respect to the Lebesgue measure, \ie if $Z \subset I$ is Borel and $\lno Z \rno = 0$ then $\mu(Z) = 0$.
By hypothesis, $\lno \leftindex^\times F(Z) \rno = 0$, \ie given $\veps > 0$ there exists a sequence $\la \hat{J}_k \ra_k$ of cells in $\R$ such that $\leftindex^\times F(Z) \subset \cup_{k=1}^\infty \hat{J}_k$ and $\sum_{k=1}^\infty \lno \hat{J}_k \rno < \veps$.
For each $k$ we consider the cell $J_k = \left( \leftindex^\times F \right)^{-1}(\hat{J}_k)$ and we notice that $Z \subset \left( \leftindex^\times F \right)^{-1}\left[\leftindex^\times F(Z)\right] \subset \cup_{k=1}^\infty J_k$.
Therefore,
\begin{equation*}
\mu(Z) \leq \sum_{k=1}^\infty F(J_k) \leq \sum_{k=1}^\infty \lno \hat{J}_k \rno < \veps,
\end{equation*}
by {\bf (i)} (specifically \eqref{eq.3}).
Since $\veps$ is arbitrary, we have estbalished that $\mu(Z) = 0$.
\par 
Since $\mu$ is a {\it finite} Borel measure, it is now classical that for every $\veps > 0$ there exists $\delta > 0$ such that $\mu(B) < \veps$ whenever $B \subset I$ is Borel and $\lno B \rno < \delta$.
In particular, if $J_1,\ldots,J_N$ are pairwise non-overlapping cells contained in $I$ and $\lno \cup_{n=1}^N J_n \rno = \sum_{n=1}^N \lno J_n \rno < \delta$ then $\sum_{n=1}^N F(J_n) = \sum_{n=1}^N \mu(J_n) = \mu \left( \cup_{n=1}^N J_n \right) < \veps$.
This completes the proof that $F$ is AC.
\par 
{\bf (iv)}
That (C) implies (A) is \ref{2.9}(A).
\end{proof}

\section{Change of variable}
\label{sec.cv}

\begin{Remark}
\label{3.1}
Suppose that $I$ and $\tilde{I}$ are cells and $\phi : I \to \tilde{I}$ is either increasing or decreasing and $\phi^\times$ is $AC$.
Since $\phi$ is monotone, we infer that it is differentiable almost everywhere.
Furthermore, either $\phi^\times$ or $-\phi^\times$ is non-negative.
Therefore, $\phi$ satisfies Luzin's condition $(N)$, by \ref{2.12}.
\end{Remark}

\begin{Empty}
\label{3.2}
Let $I$ and $\tilde{I}$ be cells and $\phi : I \to \tilde{I}$ a function.
We say that $\phi$ is {\em bi-$\boldsymbol{AC}$} if it satisfies the following properties:
\begin{enumerate}
\item[(A)] $\phi$ is a bijection;
\item[(B)] $\phi^\times$ is $AC$;
\item[(C)] $\left(\phi^{-1}\right)^\times$ is $AC$.
\end{enumerate}
Since $\phi$ is $AC$, it is continuous as well and, being a bijection, it is either increasing or decreasing and is a homeomorphism.
\end{Empty}

\begin{Theorem}
\label{3.3}
Let $I$ and $\tilde{I}$ be two cells and let $\phi : I \to \tilde{I}$ be bi-$AC$. 
If $f \in KH(\tilde{I})$ then $(f \circ \phi) \cdot |\phi'| \in KH(I)$ and 
\begin{equation*}
\int_I (f \circ \phi) \cdot |\phi'| = \int_{\tilde{I}} f.
\end{equation*}
\begin{equation*}
\begin{CD}
I @>{\phi}>> \tilde{I} @>{f}>> \R
\end{CD}
\end{equation*}
\end{Theorem}

\begin{proof}
{\bf (i)} 
We abbreviate $G = \left( \leftindex^\times F \circ \phi \right)^\times$, where $F = \left(\int^\bullet f\right)^\times$.
Thus, $G$ is an additive function defined on $\calC(I)$ and $G(J) = \rmsign(\phi') \cdot F(\phi(J))$ for all $J \in \calC(I)$, where
\begin{equation*}
\rmsign(\phi') = \begin{cases}
+1 & \text{ if $\phi$ is increasing} \\
-1 & \text{ if $\phi$ is decreasing}.
\end{cases}
\end{equation*}
In the course of this proof we shall establish that $G$ is the indefinite integral of $(f \circ \phi) \cdot \phi'$.
\par 
{\bf (ii)}
Here, we shall show that $G$ is $AC_*$.
First, we notice that $\phi$ is uniformly continuous on $I$, therefore, there exists $\boldeta : \ts ]0,\infty[ \ts\to\ts ]0,\infty[\ts$ with the following property.
For every $\delta > 0$ and every $x,y \in I$ if $|x-y| \leq \boldeta(\delta)$ then $|\phi(x)-\phi(y)| \leq \delta$.
\par 
Let $Z \subset I$ be negligible and $\veps > 0$.
According to theorem \ref{2.12}, $\phi$ satisfies Luzin's condition $(N)$, therefore, $\phi(Z)$ is negligible.
As $F$ is $AC_*$, by \ref{2.7}, there exists a gauge $\tilde{\bdelta}$ on $\phi(Z)$ such that if $\tilde{\Pi}$ is a $\tilde{\bdelta}$-fine $P$-family in $\tilde{I}$ anchored in $\phi(Z)$ then $\sum_{(\tilde{J},\tilde{x}) \in \tilde{\Pi}} |F(\tilde{J})| < \veps$.
Define a gauge $\bdelta$ on $Z$ by letting $\bdelta = \boldeta \circ \tilde{\bdelta} \circ \phi$.
Now, let $\Pi$ be a $\bdelta$-fine $P$-family in $I$ anchored in $Z$.
Define $\tilde{\Pi} = \{ (\phi(J),\phi(x)) : (J,x) \in \Pi \}$ and note that this is a $P$-family in $\tilde{I}$, since $\phi$ is increasing and continuous. 
Moreover, $\tilde{\Pi}$ is readily anchored in $\phi(Z)$ and one happily checks that it is $\tilde{\bdelta}$-fine as well, according to the definition of $\boldeta$.
Thus, $\sum_{(J,x) \in \Pi} |G(J)| = \sum_{(J,x) \in \Pi} |F[\phi(J)]| = \sum_{(\tilde{J},\tilde{x}) \in \tilde{\Pi}} |F(\tilde{J})| < \veps$.
This completes the proof that $G$ is $AC_*$.
\par 
{\bf (iii)}
Here, we shall observe that for almost every $x \in I$, $\phi$ is differentiable at $x$ and $\leftindex^\times F$ is differentiable at $\phi(x)$ with $\left( \leftindex^\times F \right)'(\phi(x)) = f(\phi(x))$, in particular, $\leftindex^\times G$ is differentiable at $x$ and $\left( \leftindex^\times G \right)'(x) = f(\phi(x)) \cdot \phi'(x)$.
Since $\phi^\times$ is $AC$, one infers that the set $Z = I \cap \{ x : \phi \text{ is not differentiable at } x\}$ is negligible (recall \ref{3.1}).
Moreover, it follows from \ref{2.7} (see \ref{2.8}(1)) that the set $\tilde{Z} = \tilde{I} \cap \{ \tilde{x} : \leftindex^\times F \text{ is not differentiable at $x$ or it is differentiable at $x$ but } \left( \leftindex^\times F \right)'(x) \neq f(x) \}$ is negligible as well.
Accordingly, $\phi^{-1}(\tilde{Z})$ is negligible, since $\phi^{-1}$ satisfies Luzin's $(N)$ condition, by \ref{3.1}.
If $x \not \in Z \cup \phi^{-1}(\tilde{Z})$ then it satisfies the claimed properties.
\par 
{\bf (iv)}
The conclusion of the theorem is a consequence of the conjunction of {\bf (ii)}, {\bf (iii)}, and \ref{2.7}:
\begin{equation*}
\begin{aligned}
\int_{\tilde{I}} f & = F(\tilde{I}) && \text{(by definition of $F$)} \\
& = F(\phi(I)) && \text{(since $\phi$ is a bijection)} \\
& = \rmsign(\phi') \cdot G(I) && \text{(by {\bf (i)})} \\
& = \rmsign(\phi') \cdot \int_I (f \circ \phi) \cdot \phi' && \text{(by \ref{2.7}, {\bf (ii)}, and {\bf (iii)}).}
\end{aligned}
\end{equation*}
\end{proof}

\section{Isometries}
\label{sec.iso}

\begin{Theorem}
\label{4.1}
If $\phi : I \to \tilde{I}$ is bi-$AC$ then $\{ \phi' = 0 \}$ is negligible.
\end{Theorem}

\begin{proof}
Of course, the set considered here is 
\begin{equation*}
A = \rmint(I) \cap \{ x : \phi \text{ is differentiable at $x$ and } \phi'(x) = 0 \}.
\end{equation*}
We observe that
\begin{equation}
\label{eq.5}
A \subset \phi^{-1} \big( \{ y : \phi^{-1} \text{ is not differentiable at } y \} \big).
\end{equation}
Indeed, if $x \in A$ and $y = \phi(x)$ then the limit of 
\begin{equation*}
\frac{\phi^{-1}(\phi(x+h)) - \phi^{-1}(y)}{\phi(x+h) - y} = \frac{h}{\phi(x+h) - \phi(x)}
\end{equation*}
as $h \to 0$ does not exist.
As $\phi(x+h) \to \phi(x)$ when $h \to 0$, this implies that $\phi^{-1}$ is not differentiable at $y$.
\par 
Since $\phi^{-1}$ is differentiable almost everywhere and satisfies Luzin's condition, \eqref{eq.5} implies that $\lno A \rno = 0$.
\end{proof}

\begin{Empty}
\label{4.2}
Given a cell $I$ we consider the vector space $KH(I)$ defined in \ref{2.2} and its vector subspace $W(I) = \KH(I) \cap \{ f : \{f \neq 0 \} \text{ is negligible} \}$. 
We let $\bKH(I)$ denote the quotient space $KH(I)/W(I)$.
For $\min I < x \leq \max I$ we abbreviate $I_x = [\min I,x]$.
If $f \in \fb \in \bKH(I)$ then, by virtue of \ref{2.8}(3), the quantity
\begin{equation*}
\| \fb \|_A = \sup \left\{ \left| \int_{I_x} f|_{I_x} \right| : \min I < x \leq \max I \right\}
\end{equation*}
is independent of the choice of $f \in \fb$.
One readily checks that $\fb \mapsto \|\fb\|_A$ is a semi-norm on $\bKH(I)$.
In view of \ref{2.8}(2), it is, in fact, a norm.
We classically call it the {\em Alexiewicz norm}, see \cite{ALE.48}.
\end{Empty}

The following is a proof of the implication $(B) \Rightarrow (A)$ of our main theorem.

\begin{Theorem}
\label{4.3}
Let $I$ and $\tilde{I}$ be two cells, $\sigma \in \{-1,1\}$, and $\phi : I \to \tilde{I}$ an increasing bi-AC function.
Then there exists a surjective linear isometry $\bT_\phi : \bKH(\tilde{I}) \to \bKH(I)$ such that $\bT_\phi(\fb) = \blseg \sigma \cdot (f \circ \phi) \cdot \phi' \brseg$ for all $\fb \in \bKH(\tilde{I})$, where $f \in \fb$.
\end{Theorem}

\begin{proof}
{\bf (i)}
First, we note that $\phi'$ is merely defined almost everywhere (recall \ref{3.1}) so that, what we really mean for the definition of $\bT_\phi(\fb)$ is that $\bT_\phi(\fb) = \blseg T_\phi(f) \brseg$, where
\begin{equation*}
T_\phi(f) : I \to \R : x \mapsto \begin{cases}
\sigma \cdot f(\phi(x)) \cdot \phi'(x) & \text{ if $\phi$ is differentiable at $x$} \\
0 & \text{ otherwise.}
\end{cases}
\end{equation*}
Observe that $T_\phi(f) \in KH(I)$, in view of \ref{3.3}.
\par 
{\bf (ii)} 
Here, we show that $\blseg T_\phi(f) \brseg$ does not depend upon the choice of $f \in \fb$.
Indeed, if $\tilde{I} \cap \{ f_1 \neq f_2 \}$ is negligible then so is $I \cap \{ f_1 \circ \phi \neq f_2 \circ \phi \} = \phi^{-1}(\tilde{I} \cap \{ f_1 \neq f_2 \})$, since $\phi^{-1}$ satisfies Luzin's condition $(N)$ (recall \ref{3.1}).
Thus, $\bT_\phi$ is well-defined.
\par 
{\bf (iii)}
One routinely checks that $\bT_\phi$ is linear.
\par 
{\bf (iv)}
Here, we show that $\bT_\phi$ is an isometry.
Let $f \in \fb \in \bKH(\tilde{I})$ and $\min I < x \leq \max I$.
Applying \ref{3.3} to $I_x$, $\tilde{I}_{\phi(x)}$, $\phi|_{I_x}$, and $f|_{\tilde{I}_{\phi(x)}}$ (notice carefully that $\phi(I_x) = \tilde{I}_{\phi(x)}$, since $\phi$ is increasing) we infer that 
\begin{equation*}
\int_{\tilde{I}_{\phi(x)}} f|_{\tilde{I}_{\phi(x)}} = \int_{I_x} \left( f|_{\tilde{I}_{\phi(x)}} \circ \phi|_{I_x}\right) \cdot \left( \phi|_{I_x}\right)' = \sigma \int_{I_x} T_\phi(f)|_{I_x}.
\end{equation*}
Accordingly, 
\begin{equation*}
\left| \int_{\tilde{I}_{\phi(x)}} f|_{\tilde{I}_{\phi(x)}} \right| = \left| \int_{I_x} T_\phi(f)|_{I_x} \right|
\end{equation*}
for all $\min I < x \leq \max I$.
Since $\phi$ is a bijection, we conclude that 
\begin{multline*}
\left\|\bT_\phi(\fb)\right\|_A
= \sup \left\{ \left| \int_{I_x} T_\phi(f)|_{I_x} \right| : \min I < x \leq \max I \right\} \\
= \sup \left\{ \left| \int_{\tilde{I}_{\phi(x)}} f|_{\tilde{I}_{\phi(x)}} \right| : \min I < x \leq \max I \right\} \\
= \sup \left\{ \left| \int_{\tilde{I}_{\tilde{x}}} f|_{\tilde{I}_{\tilde{x}}} \right| : \min \tilde{I} < \tilde{x} \leq \max \tilde{I} \right\}
= \|\fb\|_A.
\end{multline*}
\par 
{\bf (v)}
It remains to show that $\bT_\phi$ is surjective.
To this end, we first notice that the above construction applies to defining $T_{\phi^{-1}} : KH(I) \to KH(\tilde{I})$.
Given $g \in \bg \in \bKH(I)$, we note the following identities almost everywhere:
\begin{equation*}
T_{\phi^{-1}}(g) = \sigma \cdot (g \circ \phi^{-1}) \cdot (\phi^{-1})' = \sigma \cdot (g \circ \phi^{-1}) \cdot \left( \frac{1}{\phi' \circ \phi^{-1}}\right),
\end{equation*}
where the last identity almost everywhere is a consequence of \ref{4.1}.
Thus, still up to a negligible set,
\begin{equation*}
T_\phi[T_{\phi^{-1}}(g)] = \sigma^2 \cdot \big( [(g \circ \phi^{-1})] \circ \phi \big) \cdot [(\phi')^{-1}] \circ \phi \cdot \phi' = g.
\end{equation*}
This completes the proof.
\end{proof}

\begin{Empty}[Banach-Stone theorem]
\label{4.4}
We let $I$ be a cell and $C(I)$ consist of all continuous functions $F : I \to \R$.
We also let $C_0(I)$ be its subspace consisting of those $F$ such that $F(\min I) = 0$.
Equipped with the maximum norm $\|\cdot\|_\infty$, the former is a Banach space and the latter its closed subspace.
We recall Banach's theorem \cite[Ch. XI, \S4, Théorème 3 et Remarque]{BANACH} that if $S : C(\tilde{I}) \to C(I)$ is a linear surjective isometry then there exist $\sigma \in \{-1,1\}$ and a homeomorphism $\phi : I \to \tilde{I}$ such that $S(F) = \sigma \cdot (F \circ \phi)$.
Note that, {\it ibid.}, $\sigma : I \to \R$ is continuous and $|\sigma|=1$, therefore, constant.
Our aim, here, is to establish the following variant that will be used in the proof of \ref{4.6}.
\par 
{\it If $T : C_0(\tilde{I}) \to C_0(I)$ is a linear surjective isometry then there exist $\sigma \in \{-1,1\}$ and a homeomorphism $\phi : I \to \tilde{I}$ such that 
\begin{enumerate}
\item[(A)] $T(F) = \sigma \cdot (F \circ \phi)$ for all $F \in C_0(\tilde{I})$;
\item[(B)] $\phi(\min I) = \min \tilde{I}$.
\end{enumerate}
}
\par 
{\bf (i)}
Either $T(F) \geq 0$ for all $F \geq 0$ or $T(F) \leq 0$ for all $F \geq 0$.
For if not there would exist $F \geq 0$ and $G \geq 0$ such that $T(F) > 0$ and $T(G) < 0$.
For all $t > 0$, $\|T(F + t \cdot G)\|_\infty = \|F + t \cdot G\|_\infty \geq \|F\|_\infty$ and, for $t > 0$ sufficiently small, $\|T(F + t \cdot G) \|_\infty = \|T(F) + t \cdot T(G)\|_\infty < \|T(F)\|_\infty = \|F\|_\infty$.
This contradiction proves the claim. 
\par 
{\bf (ii)}
Upon observing that if the conclusions hold for $-T$ they also hold for $T$ (replacing $\sigma$ by $-\sigma$), we may henceforth assume that $T(F) \geq 0$ whenever $F \geq 0$.
\par 
{\bf (iii)}
$T(F)(I) = F(\tilde{I})$ for all $F$.
To prove this, notice that $F^+,F^- \in C_0(\tilde{I})$, $T(F^+) \geq 0$, thus, $T(F^+)(I) = [0, \|T(F^+)\|_\infty] = [0, \|F^+\|_\infty]$, and $T(-F^-) \leq 0$, therefore, $T(-F^+) = [-\|T(-F^-)\|_\infty,0] = [-\|F^-\|_\infty,0]$.
We infer that\footnote{If $a \leq s \leq 0 \leq t \leq b$ then $a \leq s + t \leq b$.} 
\begin{equation*}
T(F)(I) = [T(F^+) + T(-F^-)](I) = [-\|F^-\|_\infty,\|F^+\|_\infty] = F(I)
\end{equation*}
and the proof of the claim is complete.
\par 
{\bf (iv)} 
We notice that $C(I) \cong C_0(I) \oplus \R$ via $F \mapsto (F - F(0)) \oplus F(0)$ which allows for a linear extension $S : C(\tilde{I}) \to C(I)$ of $T$ defined by $S(F \oplus r) = T(F) \oplus r$.
As $T$ is a bijection, it easily follows that $S$ is one, too.
\par 
{\bf (v)}
Referring to {\bf (iii)} we observe that $S(F \oplus r)(I) = [T(F) \oplus r](I) = T(F)(I) + r = F(I) + r$.
Accordingly, $\|S(F \oplus r)\|_\infty = \|F \oplus r\|_\infty$.
This proves that $S$ is a linear surjective isometry and conclusion (A) follows from Banach's theorem quoted above.
\par 
{\bf (vi)}
Abbreviate $\tilde{a} = \phi(\min I)$.
Then $0 = T(F)(\min I) = \sigma \cdot F(\tilde{a})$ for all $F \in C_0(\tilde{I})$.
This readily implies that $\tilde{a} = \min \tilde{I}$.\cqfd
\end{Empty}

\begin{Empty}[Indefinite integral, again]
\label{4.5}
Recall \ref{2.3}.
By virtue of \ref{2.8}(3), we may similarly define
\begin{equation*}
\int^\bullet_I : \bKH(I) \to C_0(I)
\end{equation*}
by means of the formula
\begin{equation*}
\left( \int^\bullet_I \right) (\fb)(x) = \int_{I_x} f,
\end{equation*}
$f \in \fb$.
One checks that this is a linear map and an isometry, by the definition of the Alexiewciz norm, \ie
\begin{equation*}
\left\| \left( \int^\bullet_I \right) (\fb) \right\|_\infty = \left\| \fb \right\|_A.
\end{equation*}
\par 
We now argue that the image of $\int^\bullet_I$ is dense in $C_0(I)$.
Indeed, the set of those functions $F : I \to \R$ of the type $F = \hat{F}|_I$ for some $\hat{F} : \R \to \R$ of class $C^1$ such that $\hat{F}(\min I)=0$ is dense.
For each such $F$, $F = \left(\int^\bullet_I \right)(\blseg \hat{F}'|_I\brseg)$, by the fundamental theorem of calculus.
\end{Empty}

The following is a proof that $(A) \Rightarrow (B)$ in our main theorem.

\begin{Theorem}
\label{4.6}
Let $\bT : \bKH(\tilde{I}) \to \bKH(I)$ be a surjective linear isometry.
Then there exists $\sigma \in \{-1,1\}$ and an increasing bi-$AC$ function $\phi : I \to \tilde{I}$ such that $\bT(\fb) = \bT_\phi(\fb)$ for all $\fb \in \bKH(\tilde{I})$, where $\bT_\phi$ is as in \ref{4.3}.
\end{Theorem}

\begin{proof}
{\bf (i)}
We shall define a surjective linear isometry $T^\bullet : C_0(\tilde{I}) \to C_0(I)$ in order that the following diagram be commutative.
\begin{equation*}
\begin{CD}
\bKH(\tilde{I}) @>{\bT}>> \bKH(I) \\
@V{\int^\bullet_{\tilde{I}}}VV @V{\int^\bullet_{I}}VV \\
C_0(\tilde{I}) @>{T^\bullet}>> C_0(I)
\end{CD}
\end{equation*}
Since $\int^\bullet_{\tilde{I}}$ is injective, there is a most obvious unique way to define $T^\bullet$ on the image of $\int^\bullet_{\tilde{I}}$ so that the diagram commutes.
In this way, $T^\bullet$ restricted to the image of $\int^\bullet_{\tilde{I}}$ is linear (because so are the other three arrows) and an isometry (because so are the other three arrows).
\par 
It remains to show that $T^\bullet$ is surjective.
Let $G \in C_0(I)$.
There is a sequence $\la \bg_n \ra_n$ in $\bKH(I)$ such that $\left\| G - \left( \int^\bullet_I \right)(\bg_n)\right\|_\infty \to 0$ as $n \to \infty$, recall \ref{4.5}.
Since $\bT$ is surjective, there is a sequence $\la \fb_n \ra_n$ in $\bKH(\tilde{I})$ such that $\bT(\fb_n) = \bg_n$ for all $n$.
Abbreviate $F_n = \left( \int^\bullet_I \right)(\fb_n)$.
Since $\bT$ and both indefinite integrals involved in the above commutative diagram are isometries, we infer that $\la F_n \ra_n$ is a Cauchy sequence in $C_0(\tilde{I})$ and, therefore, converges to some $F \in C_0(\tilde{I})$.
Thus $T^\bullet(F) = \lim_n T^\bullet(F_n)$ and we are left to note that $G = \lim_n T^\bullet(F_n)$ as well:
\begin{multline*}
\lim_n \left\| G - T^\bullet(F_n) \right\|_\infty 
= \lim_n \left\| \left( \int^\bullet_I \right)(\bg_n) - \left( T^\bullet \circ \int^\bullet_{\tilde{I}} \right)(\fb_n)\right\|_\infty \\
=\lim_n \left\| \left( \int^\bullet_I \right)(\bg_n) - \left(  \int^\bullet_{\tilde{I}}  \right)(\bT(\fb_n))\right\|_\infty = \lim_n \left\| \bg_n - \bT(\fb_n) \right\|_A = 0.
\end{multline*}
\par 
{\bf (ii)}
It follows from {\bf (i)} and the Banach-Stone theorem (recall \ref{4.4}) that there exist $\sigma \in \{-1,1\}$ and an increasing homeomorphism $\phi : I \to \tilde{I}$ such that $T^\bullet(F) = \sigma \cdot (F \circ \phi)$ for all $F \in C_0(\tilde{I})$.
Furthermore, if $F \in C_0(\tilde{I})$ is in the image of $\int^\bullet_{\tilde{I}}$ then $T^\bullet(F)$ is in the image of $\int^\bullet_I$, by the commutativity of our diagram.
In other words, if $F^\times$ is $AC_*$ and $F$ is differentiable almost everywhere then $(F \circ \phi)^\times$ is $AC_*$ (and $F \circ \phi$ is differentiable almost everywhere, though, we shall not use this here), by \ref{2.7}.
We apply this observation to $F(\tilde{x}) = \tilde{x} - \min \tilde{I}$ (which is readily the indefinite integral of the constant function equal to 1 on $\tilde{I}$).
It follows that $(\phi - \min \tilde{I})^\times$ is $AC_*$ and, hence, also $\phi^\times$ itself.
As $\phi^\times$ is non-negative (because $\phi$ is increasing), we infer from \ref{2.12} that $\phi^\times$ is $AC$.
\par 
{\bf (iii)}
We note that $(T^\bullet)^{-1}(G) = \sigma \cdot (G \circ \phi^{-1})$ for all $G \in C_0(I)$.
Next, we apply the same reasoning as in {\bf (i)} and {\bf (ii)}, this time to $\bT^{-1}$ instead of $\bT$.
We infer that there exist $\tau \in \{-1,1\}$ and an increasing homeomorphism $\psi : \tilde{I} \to I$ such that $\psi^\times$ is $AC$ and such that $(T^{-1})^\bullet(G) = \tau \cdot (G \circ \psi)$ for all $G \in C_0(I)$.
If $G = \left( \int^\bullet_I \right)(\bg)$ for some $\bg \in \bKH(I)$ then
\begin{equation*}
(T^{-1})^\bullet(G) = (T^{-1})^\bullet \left( \int^\bullet_I \right)(\bg) = \left( \int^\bullet_{\tilde{I}} \right)(\bT^{-1}(\bg))
\end{equation*}
and, therefore,
\begin{equation*}
T^\bullet [ (T^{-1})^\bullet(G) ] = \left( T^\bullet \circ \int^\bullet_{\tilde{I}} \right) (\bT^{-1}(\bg)) = \left( \int^\bullet_I \circ \bT \right) (\bT^{-1}(\bg)) = \left( \int^\bullet_I \right)(\bg) = G.
\end{equation*}
Applying $(T^\bullet)^{-1}$ on both sides of the previous equation, we see that $(T^{-1})^\bullet = (T^\bullet)^{-1}$ on the image of $\int^\bullet_I$.
Choosing $G(x) = x - \min I$, $x \in I$, we infer that $\tau \cdot (\psi(\tilde{x}) - \min I) = \sigma \cdot (\phi^{-1}(\tilde{x}) - \min I)$ for all $\tilde{x} \in \tilde{I}$.
Letting $\tilde{x} = \max \tilde{I}$ we obtain that $\sigma = \tau$.
Accordingly, $\psi = \phi^{-1}$, therefore, $(\phi^{-1})^\times$ is $AC$, which concludes the proof that $\phi$ is bi-$AC$.
\par 
{\bf (iv)}
Finally, we ought to show that $\bT = \bT_\phi$.
Let $\fb \in \bKH(\tilde{I})$ and $F = \left( \int^\bullet_{\tilde{I}} \right)(\fb)$.
Then $T^\bullet(F) = \sigma \cdot (F \circ \phi)$, by {\bf (ii)}, $F' = f$ almost everywhere, by \ref{2.7} (recall \ref{2.8}(1)), and reasoning as in the proof of \ref{3.3}{\bf (iii)} we have $(F \circ \phi)' = (f \circ \phi) \cdot \phi'$ almost everywhere (since $\phi$ is bi-$AC$).
It then follows from \ref{2.7} that
\begin{equation*}
T^\bullet(F) =  \left(\int^\bullet_I \right)(\sigma \cdot \blseg (F \circ \phi)' \brseg) =  \left(\int^\bullet_I \right) ( \blseg \sigma \cdot (f \circ \phi) \cdot \phi' \brseg).
\end{equation*}
Since $T^\bullet(F) = \left( \int^\bullet_I \right) (\bT(\fb))$ as well, we conclude from the injectivity of $\int^\bullet_I$ that $\bT(\fb) = \blseg \sigma \cdot (f \circ \phi) \cdot \phi' \brseg$.
The proof is complete.
\end{proof}


\bibliographystyle{amsplain}
\bibliography{/Users/thierry/Documents/LaTeX/Bibliography/thdp}




\end{document}